\documentclass[11pt]{amsart}
\usepackage{graphicx} 
\usepackage[utf8]{inputenc}
\usepackage{amsmath,amssymb}
\usepackage{wrapfig}
\usepackage{url}
\usepackage{mathtools}
\usepackage{graphicx}
\usepackage{stmaryrd}
\usepackage{amsthm}
\usepackage{xcolor}
\usepackage[colorlinks=true,linkcolor=blue,citecolor=blue]{hyperref}
\usepackage[shortlabels]{enumitem}
\usepackage{comment}
\usepackage{relsize}
\usepackage{ dsfont }
\usepackage{stmaryrd}
\usepackage{comment}
\usepackage{geometry}
\usepackage{setspace}
\newtheorem{theorem}{Theorem}
\numberwithin{theorem}{section} 
\numberwithin{equation}{section}
\newtheorem{corollary}{Corollary}
\numberwithin{corollary}{section} 
\newtheorem{prop}{Proposition}
\numberwithin{prop}{section} 
\newtheorem{exemp}{Example}
\numberwithin{exemp}{section} 

\newtheorem{lemma}{Lemma}
\numberwithin{lemma}{section}
\theoremstyle{definition}
\newtheorem{obs}{Remark}
\theoremstyle{definition}
\newtheorem{definition}{Definition}
\numberwithin{definition}{section}

\DeclareMathOperator{\Lt}{\mathnormal{\prescript{t}{}{L}}}

\DeclareMathOperator{\supp}{\text{supp}}

\DeclareMathOperator{\T}{\mathbb{T}}

\DeclareMathOperator{\R}{\mathbb{R}}
\DeclareMathOperator{\Z}{\mathbb{Z}}
\DeclareMathOperator{\N}{\mathbb{N}}
\DeclareMathOperator{\g}{\mathbb{T}^1\times\mathbb{R}}
\DeclareMathOperator{\s}{\mathcal{S}(\g)}

\DeclareFontFamily{U}{mathx}{\hyphenchar\font45}
\DeclareFontShape{U}{mathx}{m}{n}{
      <5> <6> <7> <8> <9> <10>
      <10.95> <12> <14.4> <17.28> <20.74> <24.88>
      mathx10
      }{}
\DeclareSymbolFont{mathx}{U}{mathx}{m}{n}
\DeclareFontSubstitution{U}{mathx}{m}{n}
\DeclareMathAccent{\widecheck}{0}{mathx}{"71}
\DeclareMathOperator{\ftil}{\widehat{\widehat{\mathnormal{f}\,}}\!\!\,}
\DeclareMathOperator{\util}{\widehat{\widehat{\mathnormal{u}\,}}\!\!\,}
\DeclareMathOperator{\vtil}{\widehat{\widehat{\mathnormal{v}\,}}\!\!\,}
\DeclareMathOperator{\Lvtil}{\widehat{\widehat{\mathnormal{Lv}\,}}\!\!\,}
\DeclareMathOperator{\Ltvtil}{\widehat{\widehat{\mathnormal{\Lt v}\,}}\!\!\,}
\DeclareMathOperator{\Real}{\text{Re}}
\DeclareMathOperator{\Imag}{\text{Im}}

\newcommand{\defeq}{\vcentcolon=}

\title{Schwartz regularity of differential operators on the cylinder}
\author[A. Kowacs]{Andr\'e Pedroso Kowacs}
\address{
	Universidade Federal do Paran\'{a},
	Programa de P\'os-Gradua\c c\~ao de Matem\'{a}tica,
	C.P.19096, CEP 81531-990, Curitiba, Brazil
}
\email{andrekowacs@gmail.com}
\thanks{This study was financed in part by the Coordenação de Aperfeiçoamento de Pessoal de Nível Superior - Brasil (CAPES) - Finance Code 001. The first and third authors were supported in part by CNPq (grants 316850/2021-7 and 423458/2021-3).}
	
\subjclass{Primary 35H10, 42B05; Secondary 58D25, }

\keywords{Global solvability, Global hypoellipticity, Evolution equations, Complex vector fields}

\date{\today}

\begin{document}
\begin{abstract}
    This article presents an investigation of global properties of a class of differential operators on $\T^1\times\R$. Our approach makes use of a mixed Fourier transform, incorporating both partial Fourier series on the torus and partial Fourier transform in Euclidean space. By examining the behavior of the mixed Fourier coefficients, we obtain necessary and sufficient conditions for the Schwartz global hypoellipticity of this class of differential operators, as well as conditions for the Schwartz global solvability of said operators.  
\end{abstract}
\maketitle
\begin{singlespace}
	\tableofcontents
\end{singlespace}

\section{Introduction}
Fourier analysis has proven to be a powerful tool for solving partial differential equations not only in Euclidean space, but also in various smooth manifolds, including compact Lie groups. However, despite its wide-ranging applications, the study of the partial Fourier transform on the product of $\R^n$ with compact manifolds, such as $\T^m = (\R/2\pi\Z)^m$ remains relatively unexplored (see \cite{de_vila_Silva_2018},\cite{DEAVILASILVA2022109418},\cite{KIRILOV2020102853},\cite{KIRILOV2021108806},\cite{kowacs2023fourier}).\\
This paper focuses on investigating first order partial differential operators on $\g$, where $\T^1$ stands for the $1$-dimensional torus, also known as the $1$-dimensional sphere $\mathbb{S}^1$. Following the work in \cite{kowacs2023fourier}, by analyzing the behavior the behaviour of the partial Fourier coefficients, we deduce some global properties of differential operators on this product space. The properties we study are that of Schwartz global hypoellipticity and Schwartz global solvability, here, a differential operator $L$ is Schwartz globally hypoelliptic if whenever $Lu$ is Schwartz, then $u$ is also Schwartz. Also $L$ Schwartz global solvable if whenever $f$ is in a suitable subset of the Schwartz functions, the equation $Lu=f$ admits a Schwartz solution.\\
The motivation for this study stems from the works of F. Avila and M. Cappiello in \cite{DEAVILASILVA2022109418,silva2024systemsdifferentialoperatorstimeperiodic}, and Moraes, Kirilov and Ruzhansky in \cite{KIRILOV2020102853}, which have served as significant inspiration for our research.\\
The main results we obtained are summarized by the following theorems, whose proofs are split in the following sections.
\begin{theorem}
    Let $L$ be the differential operator on $\g$ given by
    $$L=\partial_t+c(t)\partial_x+q(t),$$
    where $c,q\in C^\infty(\T^1)$ are complex-valued smooth functions on the first variable. Then $L$ is Schwartz globally hypoelliptic if and only if the set 
    $$Z_{L_0} = \left\{(k,\xi)\in\Z\times\R | k+c_0\xi-iq_0 =0\right\}$$
    is empty and $b=\Imag(c)$ does not change sign, where 
    \begin{equation}\label{defc0q0}
        c_0=\frac{1}{2\pi}\int_0^{2\pi}c(t)dt,\quad q_0=\frac{1}{2\pi}\int_0^{2\pi}q(t)dt.
    \end{equation}
\end{theorem}
\begin{theorem}
     Let $L$ be the differential operator on $\g$ given by
    $$L=\partial_t+c(t)\partial_x+q(t),$$
    where $c,q\in C^\infty(\T^1)$ are complex-valued smooth functions on the first variable. Then $L$ is Schwartz globally solvable if $c_0=0$, $q\in i\Z$ and for every $r\in\R$, the sublevel set
    $$\Omega_r=\left\{t\in\T^1|\int_0^tb(\tau)d\tau<r\right\}$$
    is connected. Also, if $b=\Imag(c)\not\equiv0$ and $c_0\neq 0$ or $q_0\not\in i\Z$, then $L$ is Schwartz globally solvable if and only if $b=\Imag(c)$ does not change sign. Here, $c_0$ and $q_0$ are given by \eqref{defc0q0}.
\end{theorem}

 \section{Partial Fourier Analysis}
 First recall some of the definitions on the mixed partial Fourier analysis developed in \cite{kowacs2023fourier}, which will play a crucial role in proving results concerning the global properties of differential operators on $\g$.\\
 We denote by $\T^1\times\R$ the product of $\T^1=\mathbb{S}^1 = \R/2\pi\Z$, the $1$-dimensional torus and the $1$-dimensional Euclidean space $\R$. We also denote by $\mathcal{S}(\R)$ the Schwartz space on $\R$, and by $\mathcal{S}'(\R)$ its topological dual, also known as the space of tempered distributions. 
\begin{definition}\label{defipartialfourier}
    Let $f\in L^1(\g)$. We define its partial Fourier transforms by
    $$\mathcal{F}_{\R}(f)(t,\xi)=\widehat{f}(t,\xi)\defeq \int_{\R}f(t,x)e^{-ix\xi}dx,$$
    for almost every $(t,\xi)\in\g$, and
     $$\mathcal{F}_{\T^1}(f)(k,x)=\widehat{f}(k,x)\defeq \frac{1}{2\pi}\int_{0}^{2\pi}f(t,x)e^{-ik t}dt,
     $$
     for every $k\in\Z$ and almost every $x\in\R$.
     We also define its mixed partial Fourier transform by
     \begin{align*}
         \mathcal{F}_{\T^1}(\mathcal{F}_{\R}(f))(k,\xi) = \ftil(k,\xi)&\defeq\frac{1}{2\pi}\int_{0}^{2\pi}\int_{\R}f(t,x)e^{-i(k t+x\xi)}dxdt,
     \end{align*}
     for every $k\in\Z$ and almost every $x\in\R$.
     Note that by Fubini's Theorem, these are all well defined.
\end{definition}
\noindent For $g\in L^1(\g)$, we also define its inverse partial Fourier transform by
$$\mathcal{F}^{-1}_{\R}(g)(t,x) = \widecheck{g}(t,x)\defeq  \int_{\R} g(t,\xi)e^{ix \xi}d\xi,$$
and for a suitable sequence of Schwartz functions 
 $\{{g}(k,\cdot)\in C^\infty(\R)\}_{k\in\Z}$ we define its inverse partial Fourier transform by
$$\mathcal{F}^{-1}_{\T^1}(g)(t,x) = \widecheck{g}(t,x)\defeq  \sum_{k\in\Z}g(k,x)e^{ikt},$$
for almost every $(t,x)\in\T^1\times\R$.
\begin{definition}
    Let $f\in C^\infty{(\g)}$. For each $\alpha,\beta,\gamma\in\N_0$, define
    $$\|f\|_{\alpha,\beta,\gamma}\defeq\sup_{\substack{t\in\T^1\\ x\in\R}}|x^\gamma\partial_t^\alpha\partial_x^\beta f(t,x)|.$$
    Define the space of Schwartz function on $\g$ to be given by
    $$ \mathcal{S}(\g)\defeq\left\{f\in C^\infty(\g)|\|f\|_{\alpha,\beta,\gamma}<+\infty,\,\forall \alpha,\beta,\gamma\in\N_0\right\},$$
    with the topology induced by the countable family of (semi)norms
    $$p_N(f)\defeq \sum_{\alpha+\beta+\gamma\leq N}\|f\|_{\alpha,\beta,\gamma},\,N\in\N_0.$$
    It is easy to see this endows $\mathcal{S}(\g)$ a Fréchet space structure. With this in mind, we define $\mathcal{S}'(\g)$ to be its dual space, that is, the space of all linear functionals $u:\mathcal{S}(\g)\to\mathbb{C}$ such that there exist $N=N_u>0,$ and $ C=C_u>0$ such that
    $$|\langle u,f\rangle|\defeq u(f)\leq Cp_N(f),$$
    for every $ f\in \mathcal{S}(\g).$
\end{definition}
\begin{obs}\label{remarkmoderategrowth}
    Note that every $u:\T^1\times\R\to\mathbb{C}$ which satisfies
    $$|u(t,x)|\leq C(1+|x|^2)^{N/2}$$
    for some $C,N>0$ and every $(t,x)\in\g$, induces a complex measure that corresponds to an element of $\mathcal{S}'(\g)$.
\end{obs}
    
Next we recall some of the basic results obtained in \cite{kowacs2023fourier}, adjusted to the one-dimensional setting.
\begin{prop}\label{proppartialdecayrkx}
    A function $f\in C^\infty(\g)\cap L^1(\g)$ is in $\mathcal{S}(\g)$ if and only if for every $\beta,\in\N_0$, $N>0$, there exists $C_{\beta,N},R>0$, such that
    \begin{equation}{\label{ineqpartialtorus}}
        |\partial_x^\beta\widehat{f}(k,x)|\leq C_N(1+|k|^2)^{-N/2}(1+|x|^2)^{-N/2},
    \end{equation}
    for every $(k,x)\in\Z\times\R$, $|x|\geq R$. 
\end{prop}

\begin{prop}\label{proppartialdecaytxi}
    Let $f\in C^\infty(\g)$. Then $f\in\mathcal{S}(\g)\iff \mathcal{F}_{\R^{n}}(f)\in\mathcal{S}(\g)$.
\end{prop}
\noindent Notice these results imply that
$$\mathcal{F}^{-1}_{\T}(\mathcal{F}_{\T}(f))=f,$$
$$\mathcal{F}^{-1}_{\R}(\mathcal{F}_{\R}(f))=f,$$
for every $ f\in \mathcal{S}(\g).$
\begin{prop}\label{proppartialdecaykxi}
    Let $f\in C^\infty(\g)\cap L^1(\g)$. Then $f\in \mathcal{S}(\g)$ if and only if, for each $\beta\in\N_0$, $N>0$ there exist $C_{N,\beta},R>0$ such that
    \begin{equation}\label{ineqmix}
    |\partial_{\xi}^\beta\ftil(k,\xi)|\leq C_{N,\beta}(1+|k|^2)^{-N/2}(1+|\xi|^2)^{-N/2},
    \end{equation}
    for every $(k,\xi)\in\Z\times\R$, $|\xi|\geq R$.
\end{prop}

\begin{definition}
    Let $u\in \mathcal{S}'(\g)$. We define its partial Fourier transforms by
    $$\langle\mathcal{F}_{\R}(u),f\rangle = \langle u,\mathcal{F}_{\R}(f)\rangle,$$
    for every $f\in\mathcal{S}(\g)$, and for every $k\in\Z$,
    $$\mathcal{F}_{\T^1}(u)(k,\cdot)\in \mathcal{S}'(\R)$$
    is given by
    $$\langle\mathcal{F}_{\T^1}(u)(k,x),f(x)\rangle = \langle u(t,x),e^{ik\cdot t}f(x)\rangle,$$
    for every $f\in\mathcal{S}(\R)$.
    We also define its mixed partial Fourier transform by
    $$\mathcal{F}_{\T^1}(\mathcal{F}_{\R}(u))(k,\cdot)\in \mathcal{S}'(\R)$$
    $$\left\langle \mathcal{F}_{\T^1}(\mathcal{F}_{\R}(u))(k,\xi),f(\xi)\right\rangle = \langle u(t,\xi),e^{-ik\cdot t}\mathcal{F}_{\R}(f)(\xi)\rangle,$$
    for every $k\in\Z$, $f\in\mathcal{S}(\R)$.
    Also, for $u\in\mathcal{S}'(\g)$ or $\mathcal{S}'(\R)$ we define $$\langle\mathcal{F}_{\R}^{-1}u,f\rangle = \langle u,\mathcal{F}_{\R}^{-1}f\rangle.$$
    It is easy to see that $\mathcal{F}_{\R}^{-1}u\in\mathcal{S}'(\g)$ or $\mathcal{S}'(\R)$ respectively.
\end{definition}
Note the previous propositions and since $C^\infty(\T^1)\otimes\mathcal{S}(\R)\subset\mathcal{S}(\g)$, these are all well defined. Moreover, it is easy to see that this definition is consistent we Definition \ref{defipartialfourier}.
\begin{prop}\label{proppartialdecaykxdistrib}
    Let $\{\widehat{u}(k,\cdot)\in \mathcal{S}'(\R)\}_{k\in\Z}$ be a sequence of tempered distributions. Then $u\in\mathcal{S}'(\g)$ if and only if there exists $C,N>0$ such that
    \begin{equation}\label{ineqpartialtorusdistrib}
    |\langle \widehat{u}(k,x),f(x)\rangle|\leq C\Tilde{p}_N(f)(1+|k|^2)^{N/2}, 
    \end{equation}
    for every $ f\in \mathcal{S}(\R)$, where $\Tilde{p}_N(f) = \sum_{\alpha+\gamma\leq N}\|x^\gamma\partial^\alpha f\|_\infty$, for each $k\in\Z$.
\end{prop}

\begin{prop}\label{proppartialdecaytxidistrib}
     Let $u:\g\to\mathbb{C}$ such that there exists $C,N>0$ such that
    \begin{equation}\label{ineqpartialrndistrib}
    |{u}(t,\xi)|\leq C(1+|\xi|^2)^{N/2}, 
    \end{equation}
     for each $(t,\xi)\in\g$. Then there exists a unique $\widecheck{u}\in\mathcal{S}'(\g)$ such that for each $f\in \mathcal{S}(\g)$ we have that
     \begin{align*}
         \langle \widecheck{u},f\rangle &=(2\pi)^n\int_{\T^1} \int_{\R}\int_{\R}u(t,\xi)e^{ix\cdot\xi}f(t,x)d\xi dxdt\\
         &=(2\pi)^n\int_{\T^1} \int_{\R} u(t,\xi)\widehat{f}(t,-\xi)d\xi dt
      \end{align*}
     such that we define
       $\mathcal{F}^{-1}_{\R}({u}) \defeq \widecheck{u}.$
\end{prop}

\begin{prop}\label{propdecaykxidistrib}
     Let $\{{u}(k,\cdot):\R\to\mathbb{{C}}\}_{k\in\Z}$ be a sequence of functions such that $\exists N,C>0$, such that
     \begin{equation}
         |u(k,\xi)|\leq C(1+|k|^2)^{N/2}(1+|\xi|^2)^{N/2},
     \end{equation}
     for all $(k,\xi)\in\Z\times\R$. Then $\mathcal{F}^{-1}_{\R}(\mathcal{F}^{-1}_{\T^1}(u))$ is well defined and belongs to $\mathcal{S}'(\g)$. Furthermore, $\mathcal{F}^{}_{\T^1}(\mathcal{F}^{}_{\R}(\mathcal{F}^{-1}_{\R}(\mathcal{F}^{-1}_{\T^1}(u)))) = u$ in the sense of distributions.
\end{prop}

It is worth noting that as in $\mathcal{S}(\R)$, for $u\in\mathcal{S}'(\g)$, we may define: $gu$ and $\partial_t^{\alpha}\partial_x^{\beta}u$, where $g:\T^1\times\R\to\mathbb{C}$ is a suitable function (with polynomial growth) and $\alpha\in\N_0,\,\beta\in\N_0$ in the usual sense. It is also easy to see that in $\mathcal{S}'(\g)$ the following equalities also hold
\begin{align*}
    \mathcal{F}_{\T^1}(\partial_t^\alpha u)(k,x) = (ik)^\alpha \widehat{u}(k,x),\\
    \mathcal{F}_{\R}(\partial_x^\beta u)(k,x) = (i\xi)^\beta \widehat{u}(t,\xi),\\
    \mathcal{F}_{\R}((-ix)^\gamma u)(t,\xi) =  \partial_{\xi}^\gamma\widehat{u}(k,x),
\end{align*}
for any $u\in\mathcal{S}'(\g)$, $k\in\Z$, and $\alpha,\beta,\gamma\in\N_0$.

\section{Global hypoellipticity}
We now present the necessary and sufficient conditions for global hypoellipticity for a class of first order differential operators on $\g$. In the first two sections, we recall the results in obtained in \cite{kowacs2023fourier} for the constant coeffcients and variable real-valued coefficients case, and present their proofs for the sake of completion. We then present our new results for the variable complex-valued coefficients case.

\begin{definition}
    Let $L:\mathcal{S}'(\g)\to\mathcal{S}'(\g)$ be a differential operator. We say $L$ is Schwartz globally hypoelliptic (SGH) if whenever $Lu=f\in\s$ for some $u\in\mathcal{S}'(\g)$, this implies $u\in\s$.
\end{definition}
\subsection{Constant coefficients case}
\begin{theorem}\label{teoghcte}
    Let $L$ be a first order constant coefficients linear differential operator on $\T^1\times\R$, so that
    $$L = c_1\partial_t+c_2\partial_x+c_3,$$
for some $c_1,c_2,c_3\in\mathbb{C}$, $c_1\neq0$ or $c_2\neq 0$. Then $L$ is Schwartz globally hypoelliptic if and only if
$$Z\defeq \{(k,\xi)\in\Z\times\R | c_1k+c_2\xi-ic_3=0\}=\emptyset.$$
\end{theorem}
\begin{proof}
    First suppose there exists $(k_0,\xi_0)\in Z$. Then the function 
    $$v(t,x) = \frac{1}{2\pi}e^{i(k_0t+\xi_0x)}$$
    satisfies $v\in\mathcal{S}'(\T^1\times\R)\backslash\mathcal{S}(\T^1\times\R)$ and $Lv=0\in\mathcal{S}(\T^1\times\R)$, so that $L$ is not SGH.

    Now suppose $Z=\emptyset$. First, note that for any $c_0\neq 0$ $L$ is SGH if and only if $\frac{1}{c_0}L$ is SGH, and $Z=\emptyset\iff\{(k,\xi)\in\Z\times\R| \frac{c_1}{c_0}k+c\frac{c_2}{c_0}\xi-i\frac{c_3}{c_0}=0\}=\emptyset$. We  split the proof in two cases: $c_1\neq0$ and $c_2\neq0$.\\
    If $c_1\neq 0$, then by the argument above, we may suppose, without loss of generality, that 
    $$L = \partial_t+c\partial_x+q,$$
    where $c=a+ib,q\in\mathbb{C}$. If we suppose $Lu=f\in\mathcal{S}(\T^1\times\R)$, then taking the mixed partial Fourier transform on both sides of the equation yields
    \begin{align*}
    i(k+c\xi-iq)\util(k,\xi) = \ftil(k,\xi), 
    \end{align*}
    or equivalently,
    \begin{equation}\label{eqluf}
    i(k+a\xi+\Imag(q)+i(b\xi-\Real(q)))\util(k,\xi) = \ftil(k,\xi),
    \end{equation}
    for each $k\in\Z,\,\xi\in\R$.
    If $Z=\emptyset$ then then $k+c\xi-iq$ never vanishes and so 
    \begin{equation}\label{equtilf}
    \util(k,\xi) = \frac{\ftil(k,\xi)}   {i(k+c\xi-iq)}. 
    \end{equation}
    Also, from the equality \eqref{eqluf},one of the following must be true:
    \begin{itemize}
    \item $b\neq 0$ and $a\frac{\Real(q)}{b}+\Imag(q)\not\in\Z$,
    \item $b=0$ and $\Real(q)\neq0$ ,
    \item  $b=\Real(q)=a=0$ and $\Imag(q)\not\in\Z$.
\end{itemize}
    We claim that $\exists \epsilon_0>0$ such that 
    $$|k+c\xi-iq|\geq \epsilon_0,$$
    for all $(k,\xi) \in \Z\times\R$. Indeed, if $b=0,\,\Real(q)\neq 0$ or $b=\Real(q)=a=0$, $\Imag(q)\not\in\Z$ then clearly 
    $$|k+c\xi-iq|\geq \epsilon,$$
    for some $\epsilon>0$ (for instance, we can take $\epsilon = |\Real(q)|$ or $\epsilon = \min\{\lceil \Imag(q)\rceil-\Imag(q),\Imag(q)-\lfloor\Imag(q)\rfloor\}$ respectively).
    Similarly, if $b\neq0$ and $a\frac{\Real(q)}{b}+\Imag(q)\not\in\Z$, 
    there exists $\varepsilon>0$, such that 
    $$\left|k+a\frac{\Real(q)}{b}+\Imag(q)\right|\geq \varepsilon,$$ 
    for all $k\in\Z$. If $a=0$, then this implies $\Imag(q)\not\in\Z$ so $|k+c\xi-iq|\geq  \min\{\lceil \Imag(q)\rceil-\Imag(q),\Imag(q)-\lfloor\Imag(q)\rfloor\}>0$. If $a\neq0$, then for $k\in\Z$, we have that $$\left|k+a\xi+\Imag(q)\right|< \frac{\varepsilon}{2}$$
    implies that 
    $\left|\xi-\frac{\Real(q)}{b}\right|\geq\frac{\varepsilon}{2|a|}$, because otherwise
    $$\left|k+a\xi+\Imag(q)\right|\geq \left|k+a\frac{\Real(q)}{b}+\Imag(q)\right|-|a|\left|\xi-\frac{\Real(q)}{b}\right|\geq \varepsilon/2$$
    (since $a\frac{\Real(q)}{b}+\Imag(q)$ is not an integer). Therefore in this case $|k+c\xi-iq|\geq\frac{1}{2}(|k+a\xi+\Imag(q)|+|\xi-\Real(q)|)\geq \frac{1}{2}\min\left\{\frac{\varepsilon}{2},\frac{\varepsilon|b|}{2|a|}\right\}$.
    
    This way, equation (\ref{equtilf}), the previous inequality and Leibniz's formula imply that for every we have that $\beta,\gamma\in\N_0^n$, $$|\xi^\gamma\partial_{\xi}^\beta\util(k,\xi)|\leq K_{\gamma,\beta}\sup_{\gamma'\leq \gamma,\,\beta'\leq \beta}|\xi^{\gamma'}\partial_{\xi}^{\beta'}\ftil(k,\xi)|,$$ for some $K_{\gamma,\beta}>0$, so that $f\in\mathcal{S}(\T^1\times\R)\implies u\in\mathcal{S}(\T^1\times\R)$ and $L$ is SGH.\\
    Finally, if $c_1=0$, then $c_2\neq 0$ and we may suppose, without loss of generality, that 
    $$L = \partial_x+q,$$
    for some $q\in\mathbb{C}$. If we suppose $Lu=f\in\mathcal{S}(\T^1\times\R)$, then taking the mixed partial Fourier transform of both sides of the equation yields
    $$(i\xi+q)\util(k,\xi) = \ftil(k,\xi)$$
    for each $k\in\Z,\,\xi\in\R$. This way also, if $Z=\emptyset$, then $\Real(q)\neq 0$ and so $|i\xi+q|\geq |\Real(q)|>0$ so that as before, $f\in\mathcal{S}(\T^1\times\R)\implies u\in\mathcal{S}(\T^1\times\R)$ and $L$ is SGH.
\end{proof}

\begin{corollary}
    Let $L$ be a differential operator in $\T^1\times\R$ given by
    $$L = \partial_t+c\partial_x+q,$$
    where $c=a+ib,q\in\mathbb{C}$. Then $L$ is globally hypoelliptic if and only if
\begin{itemize}
    \item $b\neq 0$ and $a\frac{\Real(q)}{b}+\Imag(q)\not\in\Z$ or
    \item $b=0$ and $\Real(q)\neq0$ or
    \item  $b=\Real(q)=0=a$ and $\Imag(q)\not\in\Z$.
\end{itemize}
\end{corollary}

\subsection{Variable coefficients - real case}

In this section we will recall the result from \cite{kowacs2023fourier} which says that the differential operator
$$L = \partial_t+a(t)\partial_x+q(t),$$
where $a,q\in C^\infty(\T^1)$ and $a$ is real valued, is Schwartz globally hypoelliptic on $\T^1\times\R$ if and only if $L_0 = \partial_t+a_0\partial_x+q_0$ is Schwartz globally hypoelliptic, where 
\begin{equation}\label{eq a_0 q_0}
 a_0=\frac{1}{2\pi}\int_0^{2\pi} a(t)dt,\,q_0=\frac{1}{2\pi}\int_0^{2\pi}q(t)dt.
\end{equation}
The idea is to exhibit an operator $\Psi:\mathcal{S}'(\T^1\times\R)\to\mathcal{S}'(\T^1\times\R)$ which preserves $\mathcal{S}(\g)$ and conjugates $L$ with $L_0$. 
First, we deal with the function $a$.
\begin{prop}
    Let $a\in C^\infty(\T^1)$, $a_0\in\mathbb{R}$ be as before. Then the operator $\Psi_a:\mathcal{S}'(\T^1\times\R)\to\mathcal{S}'(\T^1\times\R)$ given by
    \begin{align*}
            (\Psi_au)(t,x) &= \frac{1}{2\pi}\int_{\R} e^{i\xi A(t)}\widehat{u}(t,\xi)e^{ix\xi}d\xi\\
            &= \mathcal{F}_{\R}^{-1}\left\{e^{i\xi A(t)}\widehat{u}(t,\xi)\right\},
    \end{align*}
    where 
    $$A(t) = \int_0^{t}a(s)ds-a_0t,$$
    is well defined and defines an automorphism over $\mathcal{S}(\T^1\times\R)$.
\end{prop}
\begin{proof}
    As $A\in C^\infty(\T^1)$ is real valued, it is easy to see that $\Psi_a$ is well defined. Also, clearly $\Psi_a$ is invertible with $\Psi_a^{-1}$ given by
    \begin{align*}
        \Psi_a^{-1}u(t,x) &=\frac{1}{2\pi}\int_{\R} e^{-i\xi A(t)}\widehat{u}(t,\xi)e^{ix\xi}d\xi\\
            &= \mathcal{F}_{\R}^{-1}\left\{e^{-i\xi A(t)}\widehat{u}(t,\xi)\right\}.
    \end{align*}
    Now let $u\in\mathcal{S}(\T^1\times\R)$, then for each $\beta\in\N_0$:
    \begin{align*}
        |\partial_\xi^\beta \widehat{\Psi_au}(t,\xi)|&=|\partial_{\xi}^\beta\left\{e^{i\xi A(t)}\widehat{u}(t,\xi\right\}|\\
        &\leq \sum_{\beta'\leq \beta}|\xi^{\beta'}\partial_{\xi}^{\beta-\beta'}\widehat{u}(t,\xi)|,
    \end{align*}
    for some $C_{\beta'}>0$, so that by Proposition (\ref{proppartialdecaytxi}) we have that $\Psi_au\in\mathcal{S}(\T^1\times\R)$. A similar argument works for $\Psi_a^{-1}$, so that $\Psi_a$ is an automorphism of $\mathcal{S}(\T^1\times\R)$.
\end{proof}
Note that $\Psi_a$ deals with the function $a$ in the following sense.
\begin{prop}
    For $a$, $\Psi_a$ as before, we have:
    $$L_{a_0}\circ\Psi_{a} = \Psi_a\circ L_a,$$
    where 
    $$L_{a_0} = \partial_t+a_0\partial_x+q(t).$$
\end{prop}
\begin{proof}
    Note that, for $u\in\mathcal{S}'(\g)$, we have that
    $$L_{a_0}\circ\Psi_au = \partial_t\Psi_au+a_0\partial_x\Psi_au+q(t)\Psi_au.$$
   Taking the partial Fourier transform in $\R$ yields
    \begin{align}
    \partial_t\widehat{\Psi_au}(t,\xi)+i\xi &a_0\widehat{\Psi_au}(t,\xi)+q(t)\widehat{\Psi_au}(t,\xi) \notag \\
    =&\partial_t\left\{e^{i\xi A(t)}\widehat{u}(t,\xi)\right\}+i\xi a_0e^{i\xi A(t)}\widehat{u}(t,\xi)+q(t)e^{i\xi A(t)}\widehat{u}(t,\xi).\label{eqautom}
    \end{align}
    But as 
    \begin{align*}
        \partial_t\left\{e^{i\xi A(t)}\widehat{u}(t,\xi)\right\} = i\xi(a(t)-a_0)e^{i\xi A(t)}\widehat{u}(t,\xi)+e^{i\xi A(t)}\partial_t\widehat{u}(t,\xi),
    \end{align*}
     equation \eqref{eqautom} implies that
    \begin{align*}
        \widehat{L_{a_0}\circ\Psi_au}(t,\xi) &= e^{i\xi A(t)}(\partial_t\widehat{u}(t,\xi)+i\xi a(t)\widehat{u}(t,\xi)+q(t)\widehat{u}(t,\xi))\\
        &=\widehat{\Psi_a\circ L u}(t,\xi),
    \end{align*}
    from which the claim follows.
\end{proof}

Next, we deal with the function $q$.
\begin{prop}
    Let $q\in C^\infty(\T^1)$, $q_0\in\mathbb{C}$ be as before. Then the operator $\Psi_{q}:\mathcal{S}'(\g)\to\mathcal{S}'(\g)$ given by
    $$(\Psi_qu)(t,x) = e^{Q(t)}{u}(t,x),$$
    where 
    $$Q(t) = \int_{0}^{t}q(s)ds-q_0t,$$
    is well defined and an automorphism of $\mathcal{S}(\g)$.
\end{prop}
The next proposition is a simple calculation, and its proof is left to the reader.
\begin{prop}
    For $q$, $\Psi_q$ as before, we have that
    $$L_0\circ \Psi_q = \Psi_q\circ L_{a_0}.$$
\end{prop}
Finally we obtain the necessary conjugation and equivalence, as follows.
\begin{corollary}\label{CoroPsi}
    Let $\Psi_a,\Psi_q$ be as before. Then 
    $$L_0\circ\Psi= L\circ \Psi$$
    and
    $$L \circ \Psi^{-1} = \Psi^{-1}\circ L_0,$$
    where $\Psi = \Psi_a\circ\Psi_q$. Therefore $L$ is Schwartz globally hypoelliptic if and only if $L_0$ is Schwartz globally hypoelliptic.
\end{corollary}
\begin{proof}
    The equalities above are an immediate consequence from the previous propositions. The proof of the second claim is as follows: note that if $L$ is SGH and $L_0u=f\in\mathcal{S}(\g)$, then $L\circ \Psi^{-1} u = \Psi^{-1}\circ L_0u=\Psi^{-1} f\in\mathcal{S}(\g)$ so that $\Psi^{-1} u\in\mathcal{S}(\g)\implies u\in\mathcal{S}(\g)$ and so $L_0$ is SGH. On the other hand, if $L_0$ is SGH and $Lu=f\in\mathcal{S}(\g)$, then $L_0\circ\Psi u = \Psi\circ L u=\Psi f\in\mathcal{S}(\g)$, so that $\Psi u\in\mathcal{S}(\g)\implies u\in\mathcal{S}(\g)$ and so $L$ is SGH.
\end{proof}

\subsection{Variable Coefficients - Imaginary case}

Now consider the class of differential operators on $\g$ given by
$$L = \partial_t+c(t)\partial_x+q(t),$$
 where $c(t)=a(t)+ib(t)$ is a smooth complex valued function on $\T^1$. Notice that to study their regularity, by the previous section we may assume, without loss of generality, that $a(t)=a_0\in\mathbb{R}$, $q=q_0\in\mathbb{C}$ are constants.\\
\subsubsection{Necessary conditions}
\begin{prop}
    If the operator $L$ is Schwartz globally hypoelliptic, then $L_0$ is Schwartz globally  hypoelliptic, where
    $$L_0 = \partial_t+c_0\partial_x+q,$$
and $c_0=\frac{1}{2\pi}\int_0^{2\pi}c(t)dt$, $q\in\mathbb{C}$.
\end{prop}
\begin{proof}
    Suppose $L_0$ is not SGH. Then by Proposition \ref{propSGH} there exists $(k_0,\xi_0)\in\Z\times\R$ such that $k_0+c_0\xi_0-iq=0$, or equivalently, $c_0\xi_0-iq\in\Z$.
    Consider $v\in\mathcal{S}'(\g)\backslash\mathcal{S}(\g)$ given by
    $$v(t,x) = \exp\left\{-\int_0^{t}i\xi_0c(s)+qds\right\}e^{i\xi_0x}.$$
    Notice that $v$ is well defined by the hypothesis on $(k_0,\xi_0)$ and it satisfies
$$Lv=0\in\mathcal{S}(\g),$$
so that $L$ is not SGH.
\end{proof}

\begin{prop}
    If the operator $L$ is Schwartz globally hypoelliptic, then $b$ does not change sign.
\end{prop}
\begin{proof}
    Suppose $b$ does change sign. If $L_0$ is not SGH, then by the previous Proposition $L$ is not SGH, and we are done. Now let's assume $L_0$ is SGH.\\
    Suppose first $b_0\geq 0$ and define:
    $$G(t,s) \defeq\int_t^{t+s}b(w)dw$$
    $$0>B\defeq \min_{s,t\in[0,2\pi]}G(t,s)=\int_{t_0}^{t_0+s_0}b(w)dw$$
    After a change of variables, we may assume, without loss of generality, that $t_0,s_0,t_0+s_0\in(0,2\pi)$. Let $\varphi\in C^\infty(\T^1)$ be such that $\supp(\varphi)\subset[t_0+s_0-\delta,t_0+s_0+\delta]$, where $\delta>0$ is small enough so that $t_0+s_0-\delta ,t_0+s_0+\delta\in(0,2\pi)$, and $0\leq \varphi\leq 1$ with $\varphi(t)\equiv1$ in $[t_0+s_0-\delta/2,t_0+s_0+\delta/2]$. Also, let $\psi\in C^\infty(\R)$ be such that $\supp(\psi)\subset[0,+\infty)$ and $\psi(\xi)\equiv1 $ in $[1,+\infty)$. For each $\xi\in\R$, define $\widehat{f}(\cdot,\xi)$ to be the $2\pi$-periodic smooth extension of
    $$t\mapsto (e^{2\pi i(\xi c_0-iq)}-1)e^{B\xi}\varphi(t)e^{-i\xi a_0(t-t_0)}e^{-q(t-t_0)}\psi(\xi).$$
    Notice that, as $B<0$ $b_0\geq 0$ and $\supp(\widehat{f}(t,\cdot))\subset [0,+\infty)$, it is easy to see that for every $\beta\in\N_0$, $N>0$, there exists $K_{\beta,N}>0$ such that:
    $$|\partial_{\xi}^{\beta}\widehat{f}(t,\xi)|\leq K_{\beta,N}(1+|\xi|^2)^{-N/2},$$
    for every $(t,\xi)\in\g$. It follows from Proposition \ref{proppartialdecaytxi} that the function given by $f = \mathcal{F}_{R}^{-1}(\widehat{f})$ is in $\mathcal{S}(\g)$. Now, notice that if $Lu=f$, then taking the partial Fourier transform on $x$ we obtain
    $$\partial_t\widehat{u}(t,\xi)+i(\xi c(t)-iq)\widehat{u}(t,\xi) = \widehat{f}(t,\xi),$$
    for $(t,x)\in\g$.    
    By Lemma \ref{lemmaanulkersol} these differential equations admit unique solution given by
    $$\widehat{u}(t,\xi) = \frac{1}{e^{2\pi i(\xi c_0-iq)}-1}\int_{0}^{2\pi}\widehat{f}(t+s)e^{\int_t^{t+s}i(\xi c(w)-iq)dw}ds,$$
    which in this case can be written as
    $$\widehat{u}(t,\xi) = e^{-i\xi a_0(t-t_0)}e^{-q(t-t_0)}\psi(\xi)\int_0^{2\pi}e^{\xi(B-G(t,s))}\varphi(t+s)ds.$$
    Notice that, as $\supp(\psi)\subset [0,+\infty)$, we have that
    $$|\widehat{u}(t,\xi)|\leq 2\pi e^{2\pi|\Real(q)|},$$
    for every $(t,\xi)\in\g$, so by Proposition \ref{proppartialdecaytxidistrib} these define $u\in\mathcal{S}'(\g)$. However, note that if $\xi$ is sufficiently big, then
    \begin{align*}
        |\widehat{u}(t_0,\xi)|&=\left|\int_0^{2\pi}e^{\xi(B-G(t_0,s))}\varphi(t_0+s)ds\right|\\
        &\geq \int_{s_0-\frac{\delta}{2}}^{s_0+\frac{\delta}{2}}e^{\xi(B-G(t_0,s))}ds\\
        &\geq K\frac{1}{\sqrt{\xi}},
    \end{align*}
    for some $K>0$, by Lemma \ref{lemmaintegralineq}. Therefore, by Proposition \ref{proppartialdecaytxi}, we conclude $u\not\in\mathcal{S}(\g)$. As we have seen that $\widehat{Lu}=\widehat{f}$, we conclude $Lu=f$ and so $L$ is not SGH, as claimed.\\
    Suppose now that $b_0<0$. The proof is similar, only now we define
    $$0<\Tilde{B} = \max_{s,t\in[0,2\pi]}\int_{t-s}^tb(w)dw=\int_{t_1-s_1}^{t_1}b(w)dw.$$
    Again, after a change of variables, we may assume, without loss of generality, that $t_1,s_1,t_1-s_1\in(0,2\pi)$. Let $\tilde{\varphi}\in C^\infty(\T^1)$ be such that $\supp(\tilde{\varphi})\subset[t_1-s_1-\delta,t_1-s_1+\delta]$, where $\delta>0$ is small enough so that $t_1-s_1-\delta ,t_1-s_1+\delta\in(0,2\pi)$, and $0\leq \varphi\leq 1$ with $\varphi(t)\equiv1$ in $[t_1-s_1-\delta/2,t_1-s_1+\delta/2]$. Then, for each $\xi\in\R$, we define $\widehat{f}(\cdot,\xi)$ to be the $2\pi$-periodic smooth extension of
    $$t\mapsto(1-e^{-2\pi i(\xi c_0-iq)})e^{-\Tilde{B}\xi}\tilde{\varphi}(t)e^{-i\xi a_0(t-t_1)}e^{-q(t-t_1)}\psi(\xi)$$
    and the rest of the proof follows similarly.
\end{proof}
\subsubsection{Sufficient condition}

\begin{lemma}\label{lemmadecayenough}
    If, for every $\xi\in\R$, $\widehat{u}(\cdot,\xi)\in C^\infty(\g)$ satisfies $\widehat{Lu}(t,\xi)=\widehat{f}(t,\xi)$, where $f\in\mathcal{S}(\g)$, and there exists $K,N,R>0$ such that
    $$|\widehat{u}(t,\xi)|\leq K (1+|\xi|^2)^\frac{N}{2}\|\widehat{f}(\cdot,\xi)\|_\infty,$$
    for every $\xi\in\R$, $|\xi|\geq R$, then $\widehat{u}\in\mathcal{S}(\g)$, so that $u\in\mathcal{S}(\g).$
\end{lemma}

\begin{proof}
    First, we claim that for each $\alpha\in\N_0$, $N>0$, there exists $C_{\alpha,N}>0$ such that
    $$|\partial_t^{\alpha}\widehat{u}(t,\xi)|\leq C_{\alpha,N}(1+|\xi|^2)^{-\frac{N}{2}},$$
    for all $(t,\xi)\in\T^1\times\R$. This can be proven inductively using the equation $\widehat{Lu} = \widehat{f}$ as follows. It is clear that the inequality holds for $\alpha=0$, as $\widehat{f}\in\mathcal{S}(\g)$. Next, suppose it holds for some $\alpha=m\in\N_0$. Then 
    \begin{align*}
        \partial_t^{m+1}\widehat{u}(t,\xi) &= \partial_t^{m}(\partial_t\widehat{u}(t,\xi))\\
        &=\partial_t^{m}(\widehat{f}(t,\xi)-(i\xi c(t)-iq)\widehat{u}(t,\xi))
    \end{align*}
    so using Leibniz's formula, the compacity of $\T^1$, the decay of $\partial_t^\alpha\widehat{f}$ and the inductive hypothesis, we obtain that the decay holds for any $\alpha\in\N_0$.\\
    Next we also claim that the decay holds for $\partial_t^{\alpha}\partial_{\xi}^{\beta}\widehat{u}(t,\xi)$ for any $\beta\in\N_0$. This is again proven by induction. Note that this holds true for $\beta=0$ as we have already seen. To see why this is true for $\beta=1$, note that, as $\widehat{Lu}=f$, we have that
    \begin{align*}
        \partial_{\xi}\widehat{f}(t,\xi) = \partial_{\xi}\partial_t\widehat{u}(t,\xi)+ic(t)\widehat{u}(t,\xi)+i\xi c(t)\partial_{\xi}\widehat{u}(t,\xi)+q\partial_{\xi} \widehat{u}(t,\xi).
    \end{align*}
    Rearranging we obtain
    $$\partial_t\partial_{\xi}\widehat{u}(t,\xi)+i(\xi c(t)-iq)\partial_{\xi}\widehat{u}(t,\xi)= \partial_{\xi}\widehat{f}(t,\xi)-ic(t)\widehat{u}(t,\xi).$$
    This means that $\partial_{\xi}\widehat{u}(t,\xi)$ satisfies the same differential equation as $\widehat{u}(t,\xi)$, only now with $\partial_{\xi}\widehat{f}(t,\xi)-ic(t)\widehat{u}(t,\xi)$ instead of $\widehat{f}(t,\xi)$. As $\partial_{\xi}\widehat{f}(t,\xi)$ and $\widehat{u}(t,\xi)$ both decay faster than any polynomial in $\xi$, we may use the same argument as before to obtain that for each $\alpha\in\N_0$, there exists $C_{\alpha,1,N}>0$ such that
    $$|\partial_t^{\alpha}\partial_{\xi}\widehat{u}(t,\xi)|\leq C_{\alpha,1,N}(1+|\xi|^2)^{-N/2},$$
    for every $(t,\xi)\in\g$. Repeating this argument inductively, we obtain the claimed result. This means, by Proposition \ref{proppartialdecaytxi}, that $u=\mathcal{F}_{\R}^{-1}\widehat{u}\in\mathcal{S}(\g)$.
\end{proof}
\begin{prop}\label{propSGH}    
    Suppose $L_0$ is Schwartz globally hypoelliptic and $b\not\equiv0$. If $b$ does not change sign, then $L$ is Schwartz globally hypoelliptic.
\end{prop}
\begin{proof}
    Suppose first $b(t)\geq 0$. This implies $ b_0>0$. As $L_0$ is SGH, in the proof of Theorem \ref{teoghcte}, we have seen that there exists $\delta>0$ such that $|k+c_0\xi-iq|\geq \delta$, for each $k\in\Z,\,\xi\in\R$. Note that this implies that there exists $\epsilon>0$ such that
    \begin{align*}
        |1-e^{-2\pi i(c_0\xi-iq)}|\geq \epsilon
    \end{align*}
    and
    \begin{align*}
        |e^{2\pi i(c_0\xi-iq)}-1|\geq \epsilon,
    \end{align*}
    for every $\xi\in\R$. If $Lu=f$, taking the partial Fourier transform in $x$ again we see that $\widehat{u}$ must be the only solution to $\widehat{Lu} = \widehat{f}$ given by the two equivalent formulas in Lemma \ref{lemmaodesol}. This way, for $\xi\geq 0$, consider the solution $\widehat{u}(t,\xi)$ given by
    \begin{align*}
        \widehat{u}(t,\xi)=\frac{1}{e^{2\pi i(\xi c_0-iq)}-1}\int_{0}^{2\pi}\widehat{f}(t+s,\xi)e^{qs}e^{i\xi a_0s-\xi\int_{t}^{t+s}b(w)dw}ds,
    \end{align*}
    and for $\xi<0$, we consider that $\widehat{u}(t,\xi)$ given by the equivalent formula
        \begin{align*}
        \widehat{u}(t,\xi)=\frac{1}{1-e^{-2\pi i(\xi c_0-iq)}}\int_{0}^{2\pi}\widehat{f}(t-s,\xi)e^{-qs}e^{-i\xi a_0s+\xi\int_{t-s}^{t}b(w)dw}ds.
    \end{align*}
    Note that in both cases it is evident that
    $$|\widehat{u}(t,\xi)|\leq \epsilon^{-1}2\pi e^{2\pi|\Real(q)|}\sup_{s\in\T}|\widehat{f}(s,\xi)|,$$
    for each $t\in\T^1$, so that by Lemma \ref{lemmadecayenough} $u\in\mathcal{S}(\g)$ and clearly $Lu=f$.
    so $L$ is SGH.\\
    The case $b(t)\leq 0$ is analogous, only now we swap the equivalent formulas for $\widehat{u}$ in each case $\xi\geq0$ and $\xi<0$.
\end{proof}
Combining the necessary and sufficient conditions obtained, we conclude:
\begin{corollary}
    The operator $L$ is Schwartz globally hypoelliptic if and only if $L_0$ is Schwartz globally hypoelliptic and $b$ does not change sign.
\end{corollary}

\section{Global Solvability}

We now investigate when does the differential equation $Lu=f$ admits solution $u\in\mathcal{S}(\g)$ for each $f\in\mathcal{S}(\g)$ satisfying certain compatibility conditions. This concept is also known as global solvability.\\
Let's go back and consider the constant coefficient differential operator given by 
$$L=c_1\partial_t+c_2\partial_x+c_3,$$
where $c_1,c_2,c_3\in\mathbb{C}$.
To determine whether or not $f\in\mathcal{S}(\g)$ can be a solution to $Lu=f$, we may look at its Fourier coefficients: note that if $Lu=f$, then
$$i(c_1k+c_2\xi-ic_3)\util(k,\xi)=\ftil(k,\xi).$$
This means that whenever $(c_1k+c_2\xi-ic_3)=0$, we must have $\ftil(k,\xi)=0$ as well. This suggests the following compatibility condition and definition of global solvability.
\begin{definition}
    Let $L$ be a constant coefficient differential operator on $\g$, so that $\Lvtil(k,\xi)=p(k,\xi)\vtil(k,\xi)$, for some polynomial $p:\Z\times\R\to\mathbb{C}$, and for every $v\in\mathcal{S}'(\g)$.
    Then we say that $L$ is globally solvable if for each $f\in\mathcal{S}(\g)$ satisfying $\ftil(k,\xi)=0$ for every $(k,\xi)\in\Z\times\R$ such that $p(k,\xi)=0$, there exists $u\in\mathcal{S}(\g)$ such that $Lu=f$.
\end{definition}
\subsection{Constant Coefficients case}
\begin{theorem}\label{teoGScte}
    For any $c_1,c_2,c_3\in\mathbb{C}$, $c_1\neq0$ or $c_2\neq0$, the operator 
    $$L=c_1\partial_t+c_2\partial_x+c_3$$
    is globally solvable.
\end{theorem}
\begin{proof}
    Let $f\in\mathcal{S}(\g)$ be such that $\ftil(k,\xi)=0$ whenever $(c_1k+c_2\xi-ic_3)=0$. 
    Suppose first that $c_1\neq 0$. We may assume, as before, that $c_1=1$, $c_2=c=a+bi,\,c_3=q$. Also, first assume $b\neq0$. If $\Imag(q)+a\frac{\Real(q)}{b}\not\in\Z$, then as in the proof of Theorem \ref{teoghcte}, we have that there exists $\epsilon>0$ such that $|\Real(k+c\xi-iq)|\geq\epsilon$, for every $(k,\xi)\in\Z\times\R$. Therefore, for each $(k,\xi)\in\Z\times\R$, we can define 
    \begin{equation}\label{equtilquot}
          \util(k,\xi)
            =\frac{\ftil(k,\xi)}{i(k+c\xi-iq)}.
    \end{equation}
    If otherwise $\Imag(q)+a\frac{\Real(q)}{b}=-k_0\in\Z$,  then for every $k\neq k_0$ we can still define $\util$ as above, and at $k=k_0$, for all $\xi\neq\frac{\Real(q)}{b}$ we can write \eqref{equtilquot} as
    \begin{align*}
\util\left(k_0,\xi\right)&=\int_0^1(\partial_\xi\ftil)\left(k_0,\theta \xi+(1-\theta)\frac{\Real(q)}{b}\right)d\theta\frac{(\xi-\frac{\Real(q)}{b})}{i(k_0+c\xi-iq)}\\
    &=\int_0^1(\partial_\xi\ftil)\left(k_0,\theta \xi+(1-\theta)\frac{\Real(q)}{b}\right)d\theta\frac{\Imag(k_0+c\xi-iq)}{ib(k_0+c\xi-iq)}
    \end{align*}
    for all $\xi\in\R$. Extending $\util$ to $\xi=\frac{\Real(q)}{b}$ by continuity through the expression above, we have that $\util$ is smooth.
   In both cases, the previous equalities, the fact that $f\in\mathcal{S}(\g)$ and Proposition \ref{proppartialdecaykxi} imply that $u\in\mathcal{S}(\g)$ and clearly $Lu=f$, so that $L$ is globally solvable. If now $b=0$ and $\Real(q)\neq0$, define $\util$ as in equation \eqref{equtilquot}. As $|\Imag(k+c\xi-iq)|>0$, $\util$ is well defined. Now if $\Real(q)=0$, and $a\neq0$, for each $k\in\Z$ we may write \eqref{equtilquot} as
     \begin{align*}
    \util\left(k,\xi\right)&=\int_0^1(\partial_\xi\ftil)\left(k,\theta \xi+(1-\theta)\frac{-k-\Imag(q)}{a})\right)d\theta\frac{(\xi-\frac{-k-\Imag(q)}{a})}{i(k+c\xi-iq)}\\
    &=\int_0^1(\partial_\xi\ftil)\left(k,\theta \xi+(1-\theta)\frac{-k-\Imag(q)}{a})\right)d\theta\frac{\Real(k+c\xi-iq)}{ia(k+c\xi-iq)}
    \end{align*}
    for $\xi\neq\frac{-k-\Imag(q)}{a}$. Therefore, defining $\util$ by the expression above and extending it by continuity over $\xi=-\frac{k+\Imag(q)}{a}$, we have that $\util$ is smooth and the claim follows as before.\\
    Now, if $a=0$, then $k-iq=0$ if and only if $k=iq\in\Z$ and so we may define $u$ by
    $$\util(k,\xi)= \frac{\ftil(k,\xi)}{i(k-iq)},$$
    whenever the denominator is not zero. Otherwise, set $\util(k,\xi)=0$.\\
    Suppose now $c_1=0$. Then $c_2\neq 0$, and we may assume $c_2=1,\,c_3=q$. It is clear that if $\Real(q)\neq0$, then
    $$\util(k,\xi)=\frac{\ftil(k,\xi)}{i\xi+q}$$
     defines a smooth fast decaying function in $\xi$ for each $k\in\Z$, such that $u\in\s$ and $Lu=f$. On the other hand, if $\Real(q)=0$, then $f(k,-\Imag(q))=0$ for each $k\in \Z$ and so we can write
    $$\util(k,\xi) = \int_0^1\partial_\xi \ftil(k,\theta\xi+(1-\theta)(-\Imag(q))d\theta\frac{\Real(\xi-iq)}{i(\xi-iq)},$$
    which again, if extended by continuity over $\xi=-\Imag(q)$, clearly defines $u\in\s$ such that $Lu=f$.
    
\end{proof}

\begin{exemp}
    In contrast with the previous theorem, if the order of $L$ isn't $1$, then $L$ might not be globally solvable. For instance, consider the operator:
    $$L = a_1\partial_t+a_2\partial_x^2$$
    where $a_1,a_2\in\R,\,a_2\neq0$. Then with the notation used in the definition of global solvability, $p(k,\xi) = ic_1k-c_2\xi^2=0\iff k=\xi=0$. Note that $f$ given by
    \begin{align*}\ftil(k,\xi)=
        \begin{cases}
            \sin(\xi)e^{-\xi^2},\quad\text{ if k=0}\\
            0,\qquad\text{ otherwise.}
        \end{cases}
    \end{align*}
    is clearly in $\mathcal{S(\g)}$ by Proposition \ref{proppartialdecaykxi}, and satisfies $\ftil(0,0)=0$, however if $Lu=f$, then $\util(0,\xi)$ can't be continuous at $\xi=0$, so that $u\not\in\mathcal{S}(\g)$ and so $L$ is not globally solvable.
\end{exemp}

\begin{prop}
    Let $L$ be a constant coefficients differential operator on $\T^1\times\R$. Then $f\in\mathcal{S}(\g)$ satisfies $\ftil(k,\xi)=0$ whenever $p(k,\xi)=0$, where $\Lvtil(k,\xi) = p(k,\xi)\vtil(k,\xi)$ if and only if $f\in(\ker \Lt)^0$, where
    $$(\ker \Lt)^0=\{f\in\mathcal{S}(\g)|\langle v,f\rangle=0,\,\forall v\in \ker\Lt\}.$$
\end{prop}
\begin{proof}
    Suppose $f$ satisfies $\ftil(k,\xi)=0$ whenever $p(k,\xi)=0$. Let $v\in\ker\Lt$. Then note that:
    $$0=\Ltvtil(k,\xi) = p(-k,-\xi)\vtil(k,\xi)$$
    so $\supp(\vtil(k,\cdot))\subset\{\xi\in\R | p(-k,-\xi)=0\}$. But then, by our hypothesis on $\ftil$, $\supp(\ftil(-k,-\cdot))\subset\{\xi\in\R| p(-k,-\xi)\neq 0\}$, so: 
    $$\langle v,f\rangle = \sum_{k\in\Z}\langle \vtil(k,\xi),\ftil(-k,-\xi)\rangle = 0$$
    and thus $f\in(\ker\Lt)^0$.\\
    Conversely, if $f\in(\ker\Lt)^0$, fix $(k_0,\xi_0)\in\Z\times\R$ such that $p(-k_0,-\xi_0)=0$, and consider $v\in\mathcal{S}'(\g)$ given by:
    \begin{align*}
    \vtil(k,\xi) = 
        \begin{cases}
            \delta_{\xi_0}(\xi),\quad\text{ if }k=k_0\\
            0,\qquad\text{ otherwise}.
        \end{cases}
    \end{align*}
    Then clearly $v\in\ker\Lt$, as
    $$\Ltvtil(k_0,\xi)=p(-k_0,-\xi)\delta_{\xi_0}(\xi)=0$$
    Therefore 
    $$0=\langle v,f\rangle=\ftil(-k_0,-\xi_0).$$
    As $(k_0,\xi)\in\Z\times\R$ such that $p(-k_0,\xi_0)=0$ was arbitrary, $f$ satisfies the claimed condition.
\end{proof}
This equivalence suggests the more general definition for global solvability, which allows us to consider this concept for differential operators with variable coefficients.
\begin{definition}
    Let $L$ be a differential operator in $\g$. We say that $L$ is globally solvable, if for each $f\in (\ker\Lt)^0$ there exists $u\in\mathcal{S}(\g)$ such that $Lu=f$.
\end{definition}
\subsection{Variable Coefficients case}

\begin{prop}
    Let $L$ be the differential operator given by
    $$L=\partial_t+a(t)\partial_x+q(t)$$
    where $a,q\in C^\infty(\T^1)$, $a$ is real-valued. Then $L$ is globally solvable.
\end{prop}

\begin{proof}
The idea is to use the operator $\Psi$ defined in Corollary \ref{CoroPsi} again. For simplicity, we will assume $q(t)=q_0\in\mathbb{C}$ is constant, but the proof for the general case is similar, only now one needs to take extra care as $\Psi^{-1}\circ\Lt_0\neq\Lt\circ\Psi^{-1}$ in general. Let $f\in(\ker\Lt)^0\cap\s$ be given. We claim that $\Psi f\in (\ker L_0^t)^0$, where 
$$L_0 = \partial_t+a_0\partial_x+q_0,$$
with $a_0,q_0$ as defined in \eqref{eq a_0 q_0}. To see this, first note that $(\ker\Lt)= \Psi^{-1}(\ker\Lt_0)$, because 
$$v\in\ker\Lt_0\iff\Lt_0v=0\iff\Psi^{-1}\Lt_0 v=0\iff\Lt\Psi^{-1}v=0\iff\Psi^{-1}v\in\ker\Lt.$$
This means that if $v\in\ker\Lt_0$, then $v=\Psi u$ for some $u\in\ker\Lt$ and so:
\begin{align*}
\langle\Psi f,v\rangle&=\langle\Psi f,\Psi u\rangle = \sum_{\xi\in\Z}\left\langle e^{i\xi A(t)}\widehat{f}(t,\xi),e^{i(-\xi)A(t)}\widehat{u}(t,-\xi)\right\rangle\\
&=\sum_{\xi\in\Z}\left\langle \widehat{f}(t,\xi),\widehat{u}(t,-\xi)\right\rangle\\
&=\langle f, u\rangle  =0
\end{align*}
by definition. Therefore $\Psi f \in(\ker\Lt_0)^0$, as claimed. Since $L_0$ is globally solvable by Theorem \ref{teoGScte}, there exists $u\in\mathcal{S}(\g)$ such that $L_0u=\Psi f$ and so $f = \Psi^{-1}L_0u=L\Psi^{-1} u$, so that $L$ is globally solvable. 
\end{proof}
Next, we consider the variable complex coefficient case. First we start with a lemma which gives a useful characterization of $(\ker \Lt)^0$.
\begin{lemma}\label{lemma0thenanul}
    Let $L$ be the differential operator given by
    $$L=\partial_t+c(t)\partial_x+q,$$
    where $c(t)=a_0+ib(t)$. Then $f\in(\ker \Lt)^0$ if $f\in\mathcal{S}(\g)$ and
    \begin{enumerate}
        \item $b_0\neq0$ and $a_0\frac{\Real(q)}{b_0}+\Imag(q)\not\in\Z$, or
        \item $b_0\neq0$, $a_0\frac{\Real(q)}{b_0}+\Imag(q)\in\Z$ and $\widehat{f}\left(\cdot,\frac{\Real(q)}{b}\right)\equiv0$ or
        \item $b_0=0$ and $\Real(q)\neq 0$ or
        \item  $b_0=0=\Real(q)$ and $\widehat{f}(\cdot,\xi)\equiv0$ whenever $a_0\xi+\Imag(q)\in\Z$,
    \end{enumerate} 
    or equivalently, $\widehat{f}(\cdot,\xi)\equiv0$ whenever $c_0\xi-iq\in\Z$.
\end{lemma}
\begin{proof}
    Note that, if $1$ or $3$ hold, then $\xi c_0+iq\not\in\Z$, for every $\xi\in\R$. If $v\in\ker\Lt$, taking the partial Fourier transform with respect to $x$ of the equation $-\Lt v=0$ yields
    $$\partial_t\widehat{v}(t,\xi)+i(\xi c(t)+iq)\widehat{v}(t,\xi) =0,$$
so Lemma \ref{lemmaodesol} implies that $\widehat{v}(\cdot,\xi)\equiv0$, for every $\xi\in\R$. Therefore $v\equiv0$ and $\forall f\in\mathcal{S}(\g)$, $f\in(\ker \Lt)^0$. On the other hand, if $2$ or $4$ hold, then by the same argument, $v\in\ker\Lt$ implies that $$\supp(\widehat{v})\subset\{(t,\xi)\in\T^1\times\R | c_0(-\xi)-iq\in\Z\}.$$
    Therefore $\widehat{v}(t,\xi)\widehat{f}(t,-\xi)=0$ almost everywhere, which implies
    $$\langle v,f\rangle = \langle \widehat{v}(t,\xi),\widehat{f}(t,-\xi)\rangle=0$$
    and so $f\in(\ker \Lt)^0$.
\end{proof}

Just as before, by the previous proposition and its proof, in what follows, from now on whenever we consider the variable complex valued coefficient  differential operator given by
$$L = \partial_t+(a(t)+ib(t))\partial_x+q(t),$$
we may assume, without loss of generality, that $a(t)\equiv a_0$, $q(t)\equiv q_0$.

\begin{prop}
    Let $L$ be as above. Suppose $b\not\equiv0$ and $b$ changes sign. If $b_0=\Real(q)=0$, also suppose $a_0\neq0$ or $\Imag(q)\not\in\Z$. Then $L$ is not globally solvable.
\end{prop} 
\begin{proof}
    First, note that if $b_0\neq0$ or $\Real(q)\neq 0$ or $a_0\neq 0$ or $\Imag(q)\not\in\Z$, then 
    $$N = \{\xi\in\R|c_0\xi-iq\not\in\Z\}$$
    is a non-empty unbounded open set. Let $\zeta\in C^\infty(\R)$ such that $0\leq \zeta(\xi)\leq 1$, $\supp(\zeta)\subset N$, and $\zeta(\xi) \equiv 1$ in $N_0\subset\Bar{N_0}\subset N$, where $N_0$ is also unbounded.
     \\
    Suppose first $b_0\geq0$. Define
    $$G(t,s) \defeq\int_t^{t+s}b(w)dw,$$
    and
    $$0>B\defeq \min_{s,t\in[0,2\pi]}G(t,s)=\int_{t_0}^{t_0+s_0}b(w)dw.$$
    After a change of variables, we may assume, without loss of generality, that $t_0,s_0,t_0+s_0\in(0,2\pi)$. Let $\varphi\in C^\infty(\T^1)$ be such that $\supp(\varphi)\subset[t_0+s_0-\delta,t_0+s_0+\delta]$, where $\delta>0$ is small enough so that $t_0+s_0-\delta ,t_0+s_0+\delta\in(0,2\pi)$, and $0\leq \varphi\leq 1$ with $\varphi(t)\equiv1$ in $[t_0+s_0-\delta/2,t_0+s_0+\delta/2]$. Also, let $\psi\in C^\infty(\R)$ be such that $\supp(\psi)\subset[0,+\infty)$ and $\psi(\xi)\equiv1 $ in $[1,+\infty)$. For each $\xi\in\R$, define $\widehat{f}(\cdot,\xi)$ to be the $2\pi$-periodic smooth extension of:
    $$t\mapsto(e^{2\pi i(\xi c_0-iq)}-1)e^{B\xi}\varphi(t)e^{-i\xi a_0(t-t_0)}e^{-q(t-t_0)}\psi(\xi)\zeta(\xi).$$
    Notice that, as $B<0$, $b_0\geq 0$ and $\supp(\widehat{f}(t,\cdot))\subset [0,+\infty)$, it is easy to see that for every $\beta\in\N_0$, $N>0$, there exists $K_{\beta,N}>0$ such that
    $$|\partial_{\xi}^{\beta}\widehat{f}(t,\xi)|\leq K_{\beta,N}(1+|\xi|^2)^{-N/2},$$
    for every $(t,\xi)\in\g$. It follows from Proposition \ref{proppartialdecaytxi} that the function given by $f = \mathcal{F}_{R}^{-1}(\widehat{f})$ is in $\mathcal{S}(\g)$. Note also that by Lemma \ref{lemma0thenanul}, $f\in(\ker \Lt)^0$, since by definition $\widehat{f}(\cdot,\xi)\equiv0$ whenever $c_0\xi-iq\in\Z$. Now, notice that if $Lu=f$, then taking the partial Fourier transform on $x$ yields
    $$\partial_t\widehat{u}(t,\xi)+i(\xi c(t)-iq)\widehat{u}(t,\xi) = \widehat{f}(t,\xi).$$
    By Lemma \ref{lemmaodesol}, for each $c_0\xi-iq\not\in\Z$, these differential equations admit unique solution given by
    $$\widehat{u}(t,\xi) = \frac{1}{e^{2\pi i(\xi c_0-iq)}-1}\int_{0}^{2\pi}\widehat{f}(t+s,\xi)e^{\int_t^{t+s}i(\xi c(w)-iq)dw}ds.$$
    Therefore, if $Lu=f$, then $u$ must satisfy 
    $$\widehat{u}(t,\xi) = e^{-i\xi a_0(t-t_0)}e^{-q(t-t_0)}\psi(\xi)\int_0^{2\pi}e^{\xi(B-G(t,s))}\varphi(t+s)ds,$$
    for each $\xi\in N_0$.
    Note also that if $\xi$ is sufficiently big, then
    \begin{align*}
        |\widehat{u}(t_0,\xi)|&=\left|\int_0^{2\pi}e^{\xi(B-G(t_0,s))}\varphi(t_0+s)ds\right|\\
        &\geq \int_{s_0-\frac{\delta}{2}}^{s_0+\frac{\delta}{2}}e^{\xi(B-G(t_0,s))}ds\\
        &\geq K\frac{1}{\sqrt{\xi}},
    \end{align*}
    for some $K>0$, by Lemma \ref{lemmaintegralineq}. Since this holds for every $\xi\in N_0,\,\xi\geq1$ and $N$ is unbounded, by Proposition \ref{proppartialdecaytxi}, we conclude $u\not\in\mathcal{S}(\g)$ and so $L$ is not globally solvable, as claimed.\\
    Suppose now that $b_0<0$. The proof is similar, but now we first define
    $$0<\Tilde{B} = \max_{s,t\in[0,2\pi]}\int_{t-s}^tb(w)dw=\int_{t_1-s_1}^{t_1}b(w)dw.$$
    Then again, after a change of variables, we may assume, without loss of generality, that $t_1,s_1,t_1-s_1\in(0,2\pi)$. Let $\tilde{\varphi}\in C^\infty(\T^1)$ be such that $\supp(\tilde{\varphi})\subset[t_1-s_1-\delta,t_1-s_1+\delta]$, where $\delta>0$ is small enough so that $t_1-s_1-\delta ,t_1-s_1+\delta\in(0,2\pi)$, and $0\leq \varphi\leq 1$ with $\varphi(t)\equiv1$ in $[t_1-s_1-\delta/2,t_1-s_1+\delta/2]$. Then, for each $\xi\in\R$, we define $\widehat{f}(\cdot,\xi)$ to be the $2\pi$-periodic smooth extension of
    $$t\mapsto(1-e^{-2\pi i(\xi c_0-iq)})e^{-\Tilde{B}\xi}\tilde{\varphi}(t)e^{-i\xi a_0(t-t_1)}e^{-q(t-t_1)}\psi(\xi)\zeta(\xi)$$
    and the rest of the proof follows similarly.
\end{proof}

\begin{prop}
    If $b\not\equiv0$, $b$ does not change sign
     then $L$ is globally solvable.
\end{prop}
\begin{proof}
    First note that as $b\not\equiv0$ and $b$ does not change sign, then $b_0\neq0$.
    First assume that $a_0\frac{Re(q)}{b_0}+\Imag(q)\not\in\Z$. This means that $c_0\xi-iq\neq0$ for all $\xi\in\R$, and the proof is the same as the proof for proving $L$ is globally hypoelliptic. Now assume $a_0\frac{\Real(q)}{b_0}+\Imag(q)\in\Z$.\\
     First note that  $c_0\xi-iq\in\Z\iff \xi = \frac{\Real(q)}{b_0}$. If $b\geq 0$ the idea is to define $\widehat{u}(t,\xi)$ to be the solution given by the two equivalent equations in Lemma \ref{lemmaodesol}.
    Note that this uniquely defines $\widehat{u}(t,\xi)$ for every $\xi\neq \frac{\Real(q)}{b_0}$. In fact, we will prove it is possible to define $\widehat{u}(t,\xi)$ at $\xi=\frac{\Real(q)}{b_0}$ in such a way that it is smooth at this point also, and so in fact such that $u\in\mathcal{S}(\g)$. First consider $\xi\leq 0$. Notice that since $f\in(\ker\Lt)^0$, we have that 
    \begin{align*}
    &\int_0^{2\pi}\widehat{f}\left(t-s,\frac{\Real(q)}{b_0}\right)e^{\frac{\Real(q)}{b_0}\int_{t-s}^tb(\tau)d\tau}e^{-ia_0\frac{\Real(q)}{b_0} s}e^{-qs}ds =\\
    &\int_{\R}\int_0^{2\pi}{f}(s,x)e^{\frac{\Real(q)}{b_0}\int_{s}^tb(\tau)d\tau}e^{-ia_0\frac{\Real(q)}{b_0} (t-s)}e^{-q(t-s)}e^{-i\frac{\Real(q)}{b_0}x}dsdx=0,
    \end{align*}
    since 
    $$\mathcal{S}'(\g)\ni v_1(s,x) = e^{\frac{\Real(q)}{b_0}\int_{s}^tb(\tau)d\tau}e^{ia_0\frac{\Real(q)}{b_0}s}e^{qs}e^{-i\frac{\Real(q)}{b_0}x}$$
    is in $\ker \Lt$, for each $t\in\T^1$, as one may verify.
    Hence, we can write
    \begin{align}
        \widehat{u}(t,\xi)&=\frac{1}{1-e^{-2\pi i (\xi c_0-iq)}}\int_0^1\int_0^{2\pi} \partial_\xi H\left(t,s,\theta\xi+(1-\theta)\frac{\Real(q)}{b_0}\right)dsd\theta\left(\xi-\frac{\Real(q)}{b_0}\right),\notag\\
       \label{equ1}
    \end{align}
    where 
    $$H(t,s,\xi) = \widehat{f}(t-s,\xi)e^{\xi\int_{t-s}^tb(\tau)d\tau}e^{-ia_0\xi s}e^{-qs}.$$
    Notice then that the quotient is bounded since the denominator
    $$1-e^{-2\pi i(\xi c_0-iq)} = 1-e^{2\pi(\xi b_0-\Real(q))}(\cos(2\pi(\xi a_0+\Imag(q)))+i\sin(2\pi(\xi a_0+\Imag(q)))$$
    vanishes only when $\xi=\frac{\Real(q)}{b_0}$ and 
    \begin{equation*}
    \lim_{\xi\to\frac{\Real(q)}{b_0}}\frac{\xi-\frac{\Real(q)}{b_0}}{1-e^{-2\pi i(\xi c_0-iq)} }=\frac{1}{2\pi ic_0},
    \end{equation*}
    from which the claim follows, as $b_0\neq 0\implies c_0\neq 0$.
    Therefore formula \eqref{equ1} defines a smooth function by continuity. It is then easy to verify, similarly to the proof of Proposition \ref{propSGH}, and using the fact that $\xi\leq0$ and $b\geq 0$ we have that
    $$|\widehat{u}(t,\xi)|\leq (2\pi)^2(1+|a_0|)\|b\|_{1} e^{2\pi|\Real(q)|}(\|\widehat{f}(\cdot,\xi)\|_\infty+\|\partial_\xi\widehat{f}(\cdot,\xi)\|_\infty),$$
    for every $(t,\xi)\in\T^1\times\R$, $\xi\leq 0$. If $\xi>0$, first note that since $f\in(\ker\Lt)^0$, we have that
    \begin{align*}
    &\int_0^{2\pi}\widehat{f}\left(t+s,\frac{\Real(q)}{b_0}\right)e^{-\frac{\Real(q)}{b_0}\int_{t}^{t+s}b(\tau)d\tau}e^{ia_0\frac{\Real(q)}{b_0} s}e^{qs}ds =\\
    &\int_{\R}\int_0^{2\pi}{f}(s,x)e^{-\frac{\Real(q)}{b_0}\int_{t}^s b(\tau)d\tau}e^{ia_0\frac{\Real(q)}{b_0} (t+s)}e^{q(t+s)}e^{-i\frac{\Real(q)}{b_0}x}dsdx=0,
    \end{align*}
    since 
    $$\mathcal{S}'(\g)\ni v_2(s,x) = e^{-\frac{\Real(q)}{b_0}\int_{t}^sb(\tau)d\tau}e^{ia_0\frac{\Real(q)}{b_0}s}e^{qs}e^{-i\frac{\Real(q)}{b_0}x}$$
    is in $\ker \Lt$, for each $t\in\T^1$, as one may verify.
    This means we can write the solution to $\widehat{Lu}(t,\xi)=\widehat{f}(t,\xi)$ as:
\begin{align*}
        \widehat{u}(t,\xi)&=\frac{1}{1-e^{-2\pi i (\xi c_0-iq)}}\int_0^1\int_0^{2\pi} \partial_\xi G\left(t,s,\theta\xi+(1-\theta)\frac{\Real(q)}{b_0}\right)dsd\theta\left(\xi-\frac{\Real(q)}{b_0}\right)\\
        &=\int_0^1\int_0^{2\pi} \partial_\xi G\left(t,s,\theta\xi+(1-\theta)\frac{\Real(q)}{b_0}\right)dsd\theta\frac{\Imag(\xi c_0-iq)\frac{1}{b_0}}{1-e^{-2\pi i (\xi c_0-iq)}},
    \end{align*}
    where 
    $$G(t,s,\xi) = \widehat{f}(t+s,\xi)e^{\xi\int_{t}^{t+s}b(\tau)d\tau}e^{ia_0\xi s}e^{qs}.$$
     Notice then that the quotient is bounded by an argument analog to the one in the previous case.
    Again, it is then easy to verify, similarly to the proof of Proposition \ref{propSGH}, and using the fact that $\xi>0$, $b\leq 0$ that
    $$|\widehat{u}(t,\xi)|\leq (2\pi)^2(1+|a_0|)\|b\|_1 e^{2\pi|\Real(q)|}(\|\widehat{f}(\cdot,\xi)\|_\infty+\|\partial_\xi\widehat{f}(\cdot,\xi)\|_\infty),$$
    for every $(t,\xi)\in\T^1\times\R$, $\xi> 0$.
    Using this and the previous estimate, we conclude by Lemma \ref{lemmadecayenough} that $u\in\mathcal{S}(\g)$ and clearly $Lu=f$. The case $b\leq 0$ is analogous, only now we switch between the cases $\xi\leq0$ and $\xi>0$. Therefore we conclude that $L$ is globally solvable. 
\end{proof}

\begin{prop}\label{propconnected}
    If $c_0=0$, $q\in i\Z$ and every sublevel
    $$\Omega_r \defeq \left\{t\in\T|\int_0^tb(s)ds<r\right\}$$
    is connected, then $L$ is globally solvable. 
\end{prop}
\begin{proof}
    Suppose $f\in(\ker\Lt)^0$. Then if $u$ is a solution to $Lu=f$, it must satisfy $\widehat{Lu}(t,\xi) = \widehat{f}(t,\xi)$ for every $t\in\T^1,\,\xi\in\R$, that is
    \begin{equation}\label{eqodeconnected}
        \partial_{t}\widehat{u}(t,\xi)+i(c(t)\xi-iq)\widehat{u}(t,\xi)=\widehat{f}(t,\xi),
    \end{equation}
    for every $t\in\T^1$, $\xi\in\R$. As $c_0=0$, $q\in i\Z$, we have that $c_0\xi-iq\in\Z$ for every $(k,\xi)\in\Z\times\R$, therefore by Lemma \ref{lemmaodesol} and Lemma $\ref{lemma0thenanul}$, for each $\xi\in\R$, we have that
\begin{equation}\label{eqwidehatuconnected}
         \widehat{u}(t,\xi) = \int_{\psi(\xi)}^te^{-i\xi\int_s^t c(\tau) d\tau+q(t-s)}\widehat{f}(s,\xi) ds
    \end{equation}
    defines a solution to the differential equation \eqref{eqodeconnected}, where $\psi\in C^{\infty}(\R)$ is such that,  $\psi(\xi) = t_1$, for all $\xi\geq 1$, $\psi(\xi)= t_2$, for all $\xi\leq -1$, where
    $$\int_0^{t_1}b(\tau)d\tau = \sup_{t\in[0,2\pi]}\int_0^tb(\tau)d\tau,$$
    and
    $$\int_0^{t_2}b(\tau)d\tau = \inf_{t\in[0,2\pi]}\int_0^tb(\tau)d\tau.$$
    This defines $\widehat{u}:\T^1\times\R\to\mathbb{C}$. We claim that 
    $$ |\widehat{u}(t,\xi)|\leq 2\pi ||\widehat{f}(\cdot,\xi)||_{\infty},$$
    for every $(t,\xi)\in\g$, $|\xi|\geq 1$.\\
    To see why, first suppose $\xi\geq 1$ and let 
    $$r_t = \int_0^t b(\tau)d\tau.$$
    Note that then $t$ and $t_0$ both belong to $\tilde{\Omega}_{r_t}$, where
    $$\tilde{\Omega}_{r_t} \defeq \left\{s\in\T|\int_0^s b(\tau)d\tau\geq r_t \right\},$$
    which is connected, by Lemma \ref{lemmasublevel}. Therefore there exists an arc of circumference $\gamma_t\subset\tilde{\Omega}_{r_t}$ connecting $t_0$ to $t$. Note that 
    $$ \int_{0}^{2\pi}e^{-i\xi\int_s^t c(\tau) d\tau+q(t-s)}\widehat{f}(s,\xi) ds=0,$$
    because $f\in(\ker\Lt)^0$ and the integral above is equal to
    $$e^{qt}\int_{\R}\int_{0}^{2\pi}e^{-ix\xi}e^{-i\xi\int_s^t c(\tau) d\tau-qs}{f}(s,x)dx ds,$$
    where $(s,x)\mapsto e^{-ix\xi}e^{-i\xi\int_s^t c(\tau) d\tau-qs}\in\mathcal{S}'(\g)$ is in $\ker\Lt$, for each $t\in\T^1$, as one may verify. This means the integral
    from $t_1$ to $t$ in equation \ref{eqwidehatuconnected} is independent on the choice of path connecting these points, so we may choose $\gamma_t$ in the definition of $\widehat{u}(t,\xi)$. Hence
    \begin{align*}
        |\widehat{u}(t,\xi)|&\leq \int_{\gamma_t}e^{\left(\int_0^t b(\tau) d\tau-\int_0^s b(\tau)d\tau\right)\xi}|\widehat{f}(s,\xi)| ds\\
        &\leq 2\pi \|\widehat{f}(\cdot,\xi)\|_{\infty}.
    \end{align*}
    The case $\xi\leq -1$ is analogous.\\
    Therefore, by Lemma \ref{lemmadecayenough} we conclude $u=\mathcal{F}^{-1}_{\R}(\widehat{u})\in\mathcal{S}(\T^1\times\R)$ and $Lu=f$, so that $L$ is globally solvable.
\end{proof}

\section{Appendix}
\begin{lemma}\label{lemmaodesol}
Let $g,\theta \in C^{\infty}(\T^1)$, and set $\theta_0 = \frac{1}{2\pi}\int_0^{2\pi}\theta(t)\mathop{dt}$. If $\theta_0\not\in i\mathbb{Z}$, then the differential equation
\begin{equation}\label{ode}
\partial_t u(t)+\theta(t)u(t)=g(t),\,\qquad t\in\T^1,
\end{equation}
admits a unique solution in $C^\infty(\T^1)$ given by
\begin{equation}\label{sol-}
u(t) = \frac{1}{1-e^{-2\pi\theta_0}}\int_0^{2\pi}g(t-s)e^{-\int_{t-s}^t\theta(\tau)d\tau}ds,
\end{equation} 
or equivalently by
\begin{equation}\label{sol+}
u(t) = \frac{1}{e^{2\pi\theta_0}-1}\int_0^{2\pi}g(t+s)e^{\int_{t}^{t+s}\theta(\tau)d\tau}ds.
\end{equation}
If $\theta_0\in i\mathbb{Z}$, then equation \eqref{ode} admits infinitely many solutions given by
\begin{equation}\label{solz}
u_\lambda(t) = \lambda e^{-\int_0^t\theta(\tau)d\tau}+\int_0^t g(s)e^{-\int_s^t\theta(\tau)d\tau}ds,
\end{equation}
for every $\lambda\in\mathbb{R}$, if and only if
$$ \int_0^{2\pi}g(t)e^{\int_0^t\theta(\tau)d\tau}\mathop{dt}=0.$$
\end{lemma}
\begin{proof}
    Since the functions on the torus may be seen as $2\pi$ periodic functions on $\R$, the proof follows from simple differentiation of the formulas above and applying the periodicity. 
\end{proof}
\begin{lemma}\label{lemmaanulkersol}
Let $L$ be the differential operator given by 
$$L=\partial_t+c(t)\partial_x+q,$$
where $c\in C^\infty(\T^1)$ and $q\in\mathbb{C}$, acting on $C^\infty(\g)$. Suppose $f\in(\ker\Lt)^0$. Then the differential equations $\widehat{Lu}(t,\xi)=\widehat{f}(t,\xi)$, that is:
$$\partial_t\widehat{u}(t,\xi) = i(\xi c(t)-iq)\widehat{u}(t,\xi)=\widehat{f}(t,\xi),$$
admit solution for every $\xi\in\R$.
\end{lemma}
\begin{proof}
    Let $c_0 = \frac{1}{2\pi}\int_0^{2\pi}c(t)dt$. Then by Lemma \ref{lemmaodesol}, if $\xi c_0-iq\not\in\Z$, then the equations above admits unique solution, for every $\widehat{f}(t,\xi)$ smooth. Now suppose $\xi c_0-iq\in\Z$. Then, notice that the function $v\in \mathcal{S'}(\g)$ given by
    $$(t,x)\mapsto e^{i\int_0^t(\xi c(\tau)-iq)d\tau}e^{-ix\xi}$$
    is in $\ker\Lt$, so
    \begin{align*}
        0 = \int_0^{2\pi}\int_{\R} e^{i\int_0^t(\xi c(\tau)-iq)d\tau}e^{-ix\xi}f(t,x)dxds=\int_0^{2\pi} \widehat{f}(t,\xi)e^{i\int_0^t(\xi c(\tau)-iq)d\tau}ds,
    \end{align*}
    and so that by Lemma \ref{lemmaodesol} the equation above admits (infinitely many) solutions.
\end{proof}

\begin{lemma}\label{lemmaintegralineq}
	Let $\psi\in C^{\infty}(\T^1)$ be a smooth real function such that $\psi(s)\geq 0$ for all $s$ and $s_0\in\T^1$ is a zero of order greater than one for $\psi$, that is, $\psi(s_0)=0=\psi'(s_0)$. Then, there exists $M>0$ such that for all $\lambda>0$ sufficiently big and $\delta>0$ we have that
	$$\int_{s_0-\delta}^{s_0+\delta}e^{-\lambda\psi(s)}ds\geq \left(\int_{-\delta}^{\delta}e^{-s^2}ds\right)\lambda^{-1/2}M^{-1/2}.$$
\end{lemma}

\begin{proof}
	By Taylor's theorem, for each $s\in (s_0-\delta,s_0+\delta)$, there exists $s'\in (s_0-\delta,s_0+\delta)$, such that
	$$\psi(s) = \frac{\psi''(s')}{2}(s-s_0)^2.$$
	Let $\tilde{M}=\sup_{s\in [s_0-\delta,s_0+\delta]}\left|\frac{\psi''(s)}{2}\right|\geq0$. If $\tilde{M}=0$, then $\psi\equiv0$ and the inequality is trivial with $M=1$. Otherwise, let $M=\tilde{M}$ and then for $\lambda M>1$ we have that
	\begin{align*}
		\int_{s_0-\delta}^{s_0+\delta}e^{-\lambda\psi(s)}ds&\geq \int_{s_0-\delta}^{s_0+\delta}e^{-(\sqrt{\lambda M}(s-s_0))^2}ds\geq \frac{1}{\sqrt{\lambda M}}\left(\int_{-\delta\sqrt{\lambda M}}^{\delta\sqrt{\lambda M}}e^{-s^2}ds\right) \\
		&\geq \frac{1}{\sqrt{\lambda M}}\left(\int_{-\delta}^{\delta}e^{-s^2}ds\right). 
	\end{align*}
\end{proof}

\begin{lemma}\label{lemmasublevel}
    Let $\phi\in C^\infty(\T^1)$ be such that $\int_0^{2\pi}\phi(t)\mathop{dt} = 0$ and for every $r\in\R$, the set $\Omega_r = \left\{t\in\T^1|\int_0^t\phi(\tau)d\tau<r\right\}$ is connected. Then the set
	$$\tilde{\Omega}_r = \left\{t\in\T^1|\int_0^t\phi(\tau)d\tau\geq r\right\}=\left\{t\in\T^1|-\int_0^t\phi(\tau)d\tau\leq -r\right\}$$
 is also connected.
\end{lemma}
\begin{proof}
	This follows from $\tilde{\Omega}_r=\T^1\backslash\Omega_r$ and the general fact that any $A\subset\T^1$ is connected if and only if $\T^1\backslash A$ is connected. To see this, let $A$ be connected. Note that the claim $\T^1\backslash A$ is connected is trivially true if $A=\T^1$. Otherwise, then $\T^1\backslash A$ contains at least one point, which, without loss of generality we may assume it is $0=0+2\pi\Z$. If we consider 
 \begin{align*}
  f:(0,2\pi)&\to\T^1\backslash\{0+2\pi\Z\}\\
  x&\mapsto x+2\pi\Z,
 \end{align*}
 then $f$ is an homeomorphism. Since $A$ is connected, $I=f^{-1}(A)$ also is. But then $I$ is an interval, so $I^c = (0,2\pi)\backslash I$ is either an interval containing $(0,\epsilon)$ or $(2\pi-\epsilon,2\pi)$ for some $\epsilon>0$ or the disjoint union of two intervals containing $(0,\epsilon)\cup(2\pi-\epsilon,0)$, for some $\epsilon>0$. In both cases, since $\T^1\backslash A = f(I^c)\cup\{0+2\pi\Z\}$ it is clearly connected. Switching the roles of $A$ and $A^c$, the converse also follows.
\end{proof}

\bibliographystyle{plain} 
\bibliography{references} 
\nocite{*}
\end{document}